\documentclass[a4paper,10pt]{scrartcl}

\usepackage{libertine}

\usepackage[T1]{fontenc}
\usepackage[utf8]{inputenc}

\usepackage{geometry}
\usepackage{amsmath,amssymb} % Für mehr mathematische Symbole
\usepackage{ntheorem} % Für neue Satz-Umgebung
\usepackage{graphicx} % Zum Einbinden von Graphiken
\usepackage{enumerate}
\usepackage{enumitem}
\usepackage{tikz}

\usepackage{faktor}
\usepackage{microtype}

\newcommand{\ff}{\mathbb F}

\newcommand{\nn}{\mathbb N}

\newcommand{\qq}{\mathbb Q}
\newcommand{\rr}{\mathbb R}
\newcommand{\cc}{\mathbb C}

\newcommand{\fakt}[2]{#1\,\raise1pt\hbox{\ensuremath{/}}#2} %Definiere faktor um

\theoremheaderfont{\bfseries\sffamily}
\theorembodyfont{\rmfamily}
\theoremseparator{:}
\newtheorem*{bew}{Proof}
\newenvironment{beweis}{\begin{bew}}{\qed\end{bew}\vspace{2mm}}

\theoremheaderfont{\bfseries\sffamily}
\theorembodyfont{\rmfamily}
\theoremseparator{.}
\newcommand{\qed}{\hfill\ensuremath{\square}}

\newcounter{count}
\newtheorem{Satz}[count]{Theorem}

\newtheorem{Lemma}[count]{Lemma}

\newtheorem{Bemerkung}[count]{Remark}

    \makeatletter
    \renewcommand*{\@fnsymbol}[1]{\ensuremath{\ifcase#1\or $§$\or \sharp\or \ddagger\or
       \mathsection\or \mathparagraph\or \|\or **\or \dagger\dagger
       \or \ddagger\ddagger \else\@ctrerr\fi}}
    \makeatother

\title{\LARGE Septic equations are solvable  by 2-fold origami}
\author{ Joachim König\thanks{joachim.koenig@mathematik.uni-wuerzburg.de}, \, Dmitri Nedrenco\thanks{dmitri.nedrenco@mathematik.uni-wuerzburg.de}}

\begin{document}

\maketitle
\begin{abstract}
\noindent In this paper we prove that a generic rational equation of degree $7$ is solvable by 2-fold origami. In particular we show how to septisect an arbitrary angle.
This extends the work of \cite{AL} and \cite{Ni} on 2-fold origami.
Furthermore we give exact crease patterns for folding polynomials with Galois groups $A_7$ resp. $PSL_3\ff_2$. 
\end{abstract}

\section{Motivation}
Almost every paper about geometry can start with ``Ancient Greeks already knew how to\ldots''. Ancient Greeks knew how to trisect an angle, for example by neusis \cite[Theorem 9.3]{martin} or with a conchoid of Nicomedes \cite{bries}. But they could not trisect an angle with a ruler and compasses. After Wantzel we know that it is indeed impossible \cite{wantzel}. In the 1930s Margherita Beloch found out that one can trisect an angle by paper folding \cite{beloch},  \cite{hullbel}. Her ideas were almost forgotten and a new wave of origamists was needed to describe the power of paper folding. By the end of the second millennium it was proven that paper folding can solve arbitrary (rational) quartic polynomials, so every 2-3-tower over $\qq$ is constructible by means of paper folding. By \emph{paper folding} we mean the so called 1-fold origami; only one foldline is allowed in each folding step, cf. \cite{AL}. Generalising this, one defines $n$-fold origami by allowing $n$ fold lines (simultaneously) arising in each folding step. In 2006 Alperin and Lang developed axioms for 2-fold origami and calculated ideals describing each of the axioms: two simultaneous fold lines can be produced in every folding step (think, for instance, of folding a letter). They proved (cf. \cite[Theorem 1]{AL}), using the method of Lill \cite{hullbel}, \cite{lill}, that 
\begin{Satz}\label{n-2}
Every polynomial of degree $n$ can be solved by $(n-2)$-fold origami.
\end{Satz} So in particular you need at most 3-fold origami to solve quintics. Alperin and Lang asked whether you can do better. Nishimura showed that every quintic is solvable by means of one 2-fold axiom (AL4a6ab in the Alperin and Lang notation). He did it by interpreting this axiom geometrically (and some quite involved calculations). It is a remarkable improvement of Theorem \ref{n-2}. We try to take this game a little bit further and investigate whether one can solve every septic equation with one or more 2-fold axioms.
% 
% \footnote{In particular this solves the problem of finding a seventh root of 2 and septisection of an arbitrary angle by 2-fold.}. 
% 
We see this work as a continuation of the papers \cite{AL} and \cite{Ni}.
\begin{figure}[!ht]
\centering
\includegraphics[height=8cm]{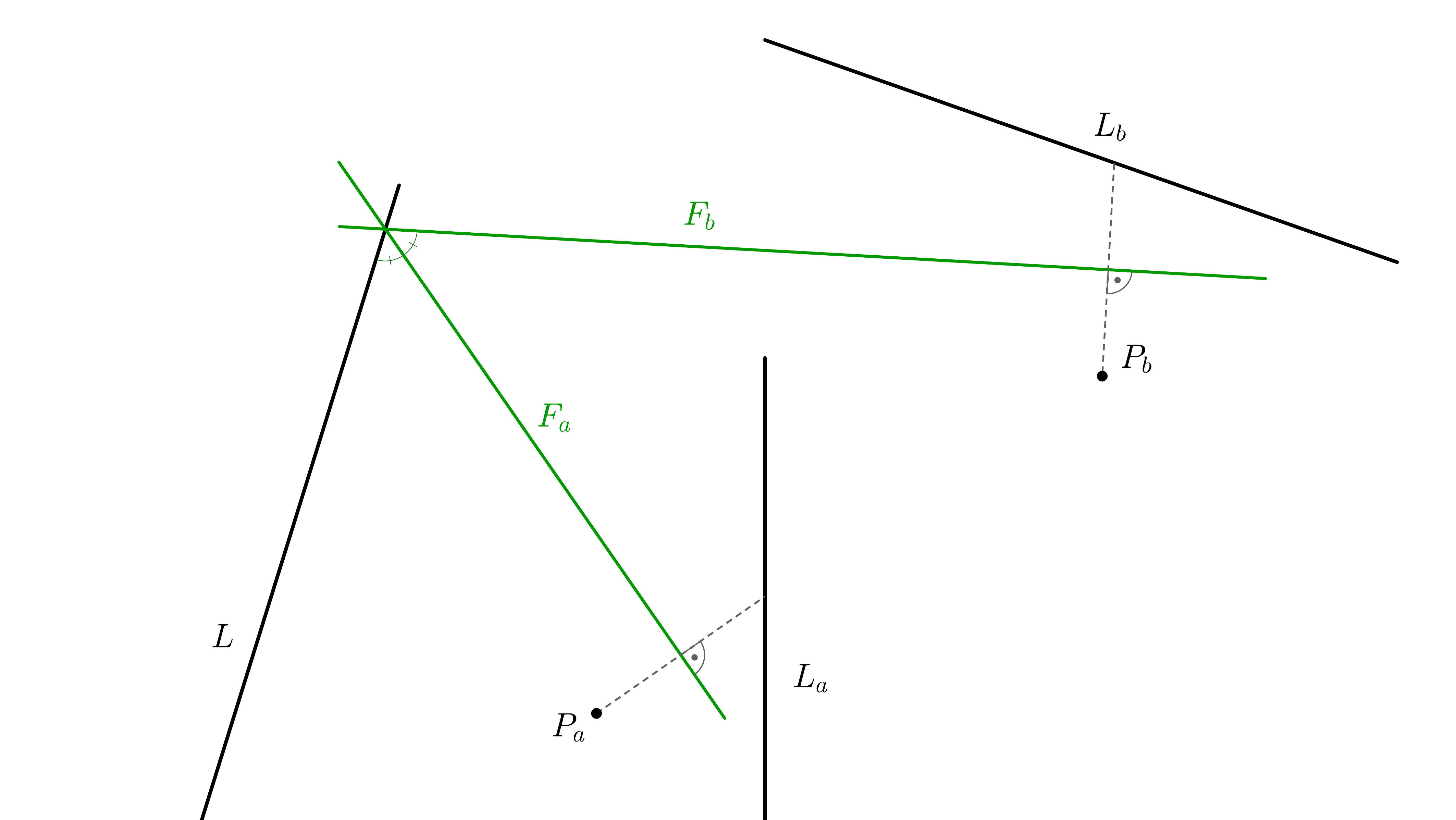}
\caption{Axiom AL4a6ab used by Nishimura. $L^{F_a}=F_b$,\; $P_a^{F_a}\in L_a$,\; $P_b^{F_b}\in L_b$. Notation of the points and lines as in \cite{Ni}.}
\end{figure}

\section{Setting}

We use the notation of a \emph{point}, \emph{line}, \emph{folded point} and \emph{folded line} as Alperin and Lang do, cf. \cite[pp.\ 4–5]{AL} for exact definitions and formulas.

The list of 2-fold axioms given by Alperin and Lang is impressive and too long in order to try every axiom. So we did some calculations in order to filter out the axioms which yield irreducible polynomials of degree $7$ for the slope of one of the fold lines: These are\footnote{One might think of ``AL'' as Alignment or Alperin-Lang.}: AL3a5b6b7a, AL3a5b6b7b, AL3a5b7ab and AL6ab8. From the geometrical point of view and after consulting \cite[7.1.1]{AL} we decided that AL6ab8 is suitable. Let us describe it.

Assume that we have already constructed four points and two lines. We seek to fold one point onto the first line and the second point onto the second line such that the third point folded by the first foldline meets the fourth point folded by the second foldline, cf. Figure 2.

\begin{figure}[!ht]
\centering
\includegraphics[height=7cm]{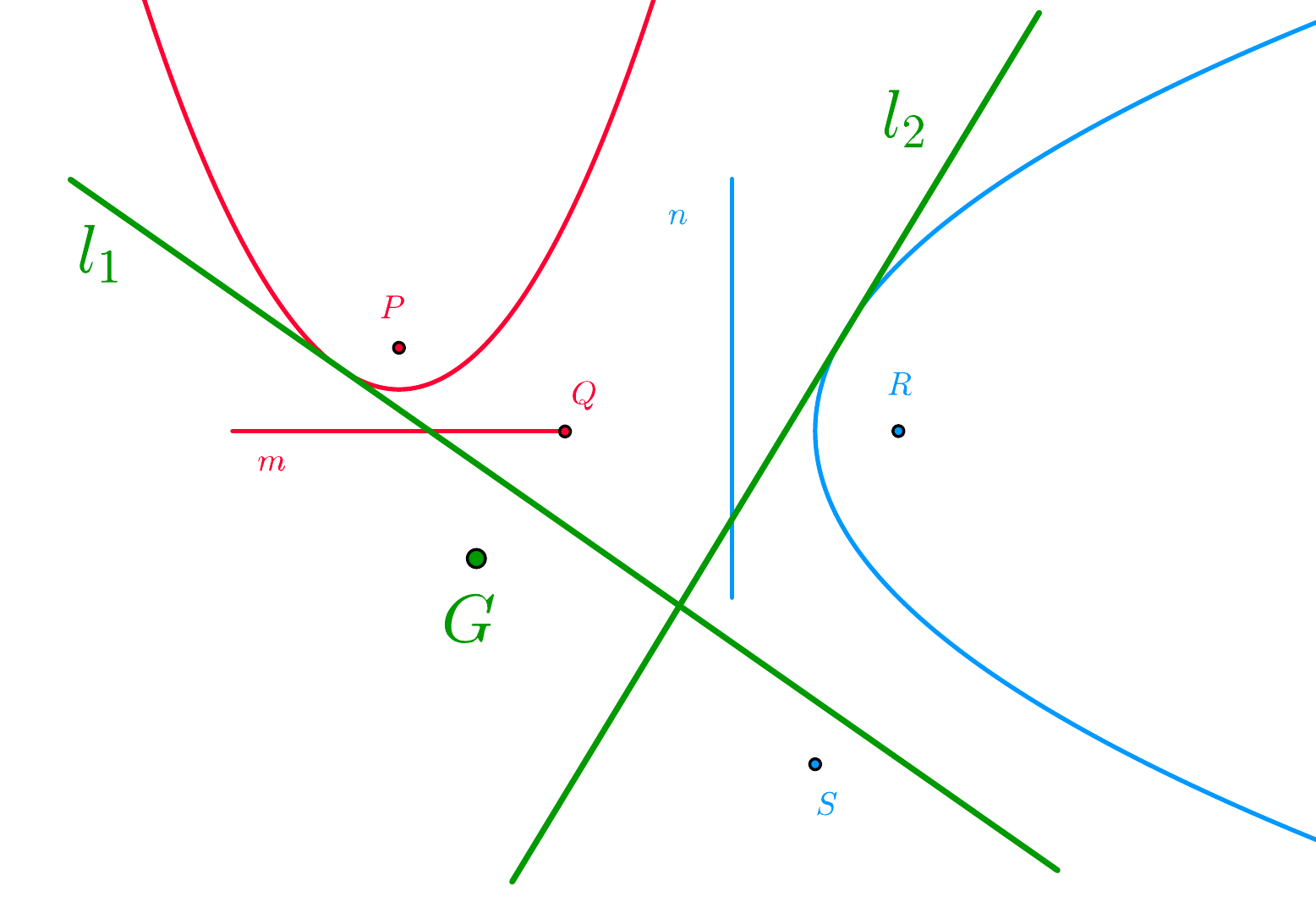}
 \caption{A representation of the axiom AL6ab8. $Q^{l_1} = G = S^{l_2}$, $P^{l_1}\in m$, $R^{l_2}\in n$.}
\end{figure}

From the basic origami theory we know that folding one point not on a given line onto this line yields as a fold line a tangent to the parabola defined by the point and line \cite[Theorem 10.3]{martin}.
We fix the line $m$ and the points $P\not\in m$ and $Q\neq P$ and we let the tangent $l$ to the parabola with focus $P$ and directrix $m$  vary over the set of all tangents to this parabola. Then the reflexion image $Q^{l}$ of $Q$ across $l$ moves along a cubic curve, cf. Figure 3.  This curve was discussed for instance in \cite[p.\ 150]{martin}, \cite[pp.\ 76]{hull}, \cite{frigerio}.

\begin{figure}[!ht]\centering
\includegraphics[height=7cm]{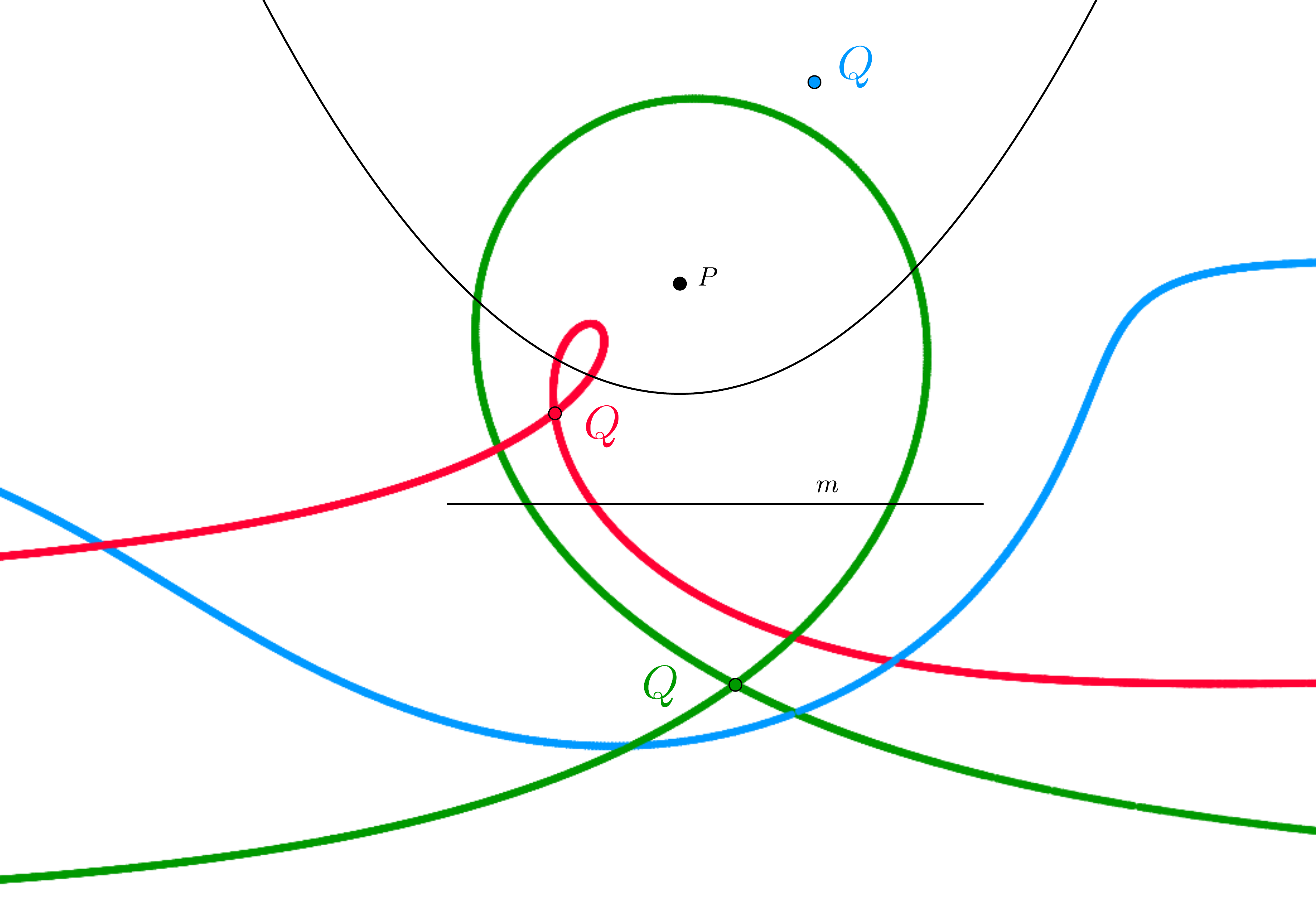}
\caption{Three different origami cubic curves for three different positions of $Q$. A point $Q$ is folded across each tangent to the parabola with directrix $m$ and focus $P$. Note that the blue curve has an isolated \emph{real} point $Q$, but obviously there are no such isolated points if you look at this curve in $\mathbb{P}_2\cc$ where it naturally lives.  }
\end{figure}

% \footnote{Called \emph{(circular nodal) origami cubic curve} in \cite{AL}.}

Let us be more specific here. Let $m =\{(x,y)\in\rr^2 \mid ax+by+1=0\}$ be a given line with origami-constructible\footnote{We will use \emph{origami-constructible} for \emph{1-fold origami-constructible}.} $a,b$. Let $P=(c,d)\in\rr^2, P\not\in m$ and $Q=(e,f)\in\rr^2$, $Q\neq P$ be two given points with origami-constructible coefficients $c,d,e,f$. It is well-known how to find the image of a point by reflexion across a line, cf. \cite[Definition 3]{AL}. The equation of an arbitrary tangent $l$ to a parabola with focus $P$ and directrix $m$  is easily calculated, too. So if we put this together, we can calculate the equation of the moving point $Q^{l}$. It is: $ax^3 + bx^2y + (-ac - ae + bd - bf + 1)x^2 + axy^2 + (-2ad - 2bc)xy + (2ace + 2adf - ae^2 - af^2 + 2bcf - 2bde - 2e)x + by^3 + (ac - ae - bd - bf + 1)y^2 + (-2acf + 2ade + 2bce + 2bdf - be^2 - bf^2 - 2f)y - ace^2 +
    acf^2 - 2adef + ae^3 + aef^2 - 2bcef + bde^2 - bdf^2 + be^2f + bf^3 + e^2 + f^2 = 0 $.
Obviously, this curve passes through $Q$, so if we change the coordinate system to $X\!:=x-e$ and $Y\!:=y-f$ and simplify a little bit the curve will have the equation:
\begin{equation}\label{cubic}
 (X^2+Y^2)(aX+bY) + t_1 X^2+t_2 XY +t_3 Y^2 = 0
\end{equation}
with $t_1=-ac + 2ae + bd + 1, t_2=-2ad + 2af - 2bc + 2be, t_3=ac - bd + 2bf + 1\in\qq(a,b,c,d,e,f)$. 
We call this curve a \emph{(circular nodal)} origami cubic curve and denote it by $\mathcal{C}\!:=\mathcal{C}(a,b,c,d,e,f)$.

The letters $a,b,c,d,e,f$ will be used throughout the paper in the sense just explained.
We pause for a second to show a converse to the necessary condition for an origami cubic curve.

In the definition of the line $m=\{(x,y)\mid ax+by+1=0\}$ there is no loss of generality assuming $a\neq 0$. Then we divide the equation~(\ref{cubic}) by $a$ and set 
$t_0:=\frac{b}{a}$ getting
\begin{equation*}
 (X^2+Y^2)(X+t_0Y) + (-c + 2e + t_0d + \frac{1}{a})X^2 + (-2d + 2f - 2t_0c + 2t_0e)XY + (c - t_0d + 2t_0f + \frac{1}{a})Y^2 = 0.
\end{equation*}
We show that the coefficients at $X^2$, $XY$, and $Y^2$ can take arbitrary values $t_1$, $t_2$, $t_3$ (for given $t_0,e,f$ and suitable choices of $a,c,d$).
To achieve this, we have to solve the system of equations
\begin{align*}
&t_0= \frac{b}{a},\\
&t_1=-ac + 2ae + t_0ad + 1,\\
&t_2=-2d + 2f - 2t_0c + 2t_0e ,\\
&t_3=ac - t_0ad + 2t_0af + 1
\end{align*}
for given $t_0,t_1,t_2,t_3,e,f$ and unknown variables $a,b,c,d$.
This system is solved by
\begin{align*}
&a:=\frac{2}{2t_0f - t_1 - t_3 + 2e},\\
&b:=t_0 a,\\
&c:=\frac{2t_0^2e - t_0 t_2 - t_1 + t_3 + 2e}{2t_0^2 + 2},\\
&d:=\frac{2t_0^2 f + t_0 t_1 - t_0 t_3 - t_2 + 2f}{2t_0^2 + 2}.
\end{align*}

We have therefore shown the following 
\begin{Lemma}
Let $t_0,t_1,t_2,t_3,e,f$ be real numbers, and let $X\!:=x-e$ and $Y\!\!:=y-f$. Then the cubic curve given by $(X^2+Y^2)(X+t_0 Y) +t_1X^2+t_2XY +t_3Y^2=0$ is an origami cubic curve $\mathcal{C}(a,b,c,d,e,f)$
with $a,b,c,d$ in the field generated by $t_0,t_1,t_2,t_3,e,f$ over $\qq$. That is, there exists a parabola (given by directrix $\{(x,y)\mid ax+by+1=0\}$ and focus $(c,d)$) such that the image
of the given point $(e,f)$ under reflection across the tangents of the parabola is exactly the given cubic.
\end{Lemma}

\begin{Bemerkung}
For certain values of $t_0,t_1,t_2,t_3,e,f$ the denominator of $a$ (and $b$) in the above solution may vanish. This, however, does not mean that there is no solution, 
but rather that the directrix of the parabola passes through the origin and therefore cannot be represented in the form $g_{a,b}:=\{(x,y)\in~\!\rr^2, ax+by+1=0\}$ as above.
Furthermore one can see that if the cubic curve is irreducible then the point $(c,d)$ in the above solution will not lie on the line $g_{a,b}$, i.e.\! they really define a parabola. 
\end{Bemerkung}
% \begin{Lemma}
%  Every cubic curve $C$ of the form $(X^2+Y^2)(X+t_0Y) + t_1X^2+t_2XY +t_3Y^2=0$ with $t_0,t_1,t_2,t_3\in\qq(a,b,c,d,e,f)$ is a cubic origami curve defined over $\qq(a,b,c,d,e,f)$. 
% \end{Lemma}
% 
% \begin{beweis}
%  As $(a,b)$ defines the line $\{(x,y)\mid ax+by+1=0\}$ there is no loss of generality assuming $a\neq 0$. Then we divide the equation~(\ref{cubic}) by $a$ and set $t_0:=\frac{b}{a}$ getting
% \begin{equation}
%  (X^2+Y^2)(X+t_0Y) + (-c + 2e + t_0d + \frac{1}{a})X^2 + (-2d + 2f - 2t_0c + 2t_0e)XY + (c - t_0d + 2t_0f + \frac{1}{a})Y^2 = 0.
% \end{equation}
% We have to show that the coefficients at $X^2$, $XY$, and $Y^2$ can take arbitrary values $t_1$, $t_2$, $t_3$ (for given $t_0,e,f$ and suitable choices of $a,c,d$).
%  The arising system of equations is solved by 
% \begin{align*}
% a:=&1 \Big/ \Big(-t_0f + \frac{1}{2}t_1 + \frac{1}{2}t_3 - e\Big),\\
%  c:=& \Big(b^2t_1^2e - b^2t_1t_2f + 2b^2t_1t_3e - 4b^2t_1e^2 - 2b^2t_1f^2 - b^2t_2t_3f + 2b^2t_2ef +\\&+ b^2t_3^2e - 4b^2t_3e^2 + 2b^2t_3f^2 + 4b^2e^3 + 4b^2ef^2 - bt_1t_2 - 4bt_1f - bt_2t_3 + 2bt_2e + 4bt_3f + 8bef -\\&- 2t_1 + 2t_3 + 4e\Big) \Big/ \Big(b^2t_1^2 + 2b^2t_1t_3 - 4b^2t_1e + b^2t_3^2 - 4b^2t_3e + 4b^2e^2 + 4b^2f^2 + 8bf + 4\Big),\\
%  d :=& \Big(2b^2t_1^2f + 2b^2t_1t_3f - 6b^2t_1ef - 2b^2t_2f^2 - 2b^2t_3ef + 4b^2e^2f + 4b^2f^3 + bt_1^2 - 2bt_1e - \\&- 4bt_2f - bt_3^2 + 2bt_3e + 8bf^2 - 2t_2 - 4f\Big) \Big/ \Big(b^2t_1^2 + 2b^2t_1t_3 - 4b^2t_1e + b^2t_3^2 - 4b^2t_3e+ \\&+ 4b^2e^2 + 4b^2f^2 +
%         8bf + 4\Big).
% \end{align*}
% This solution concludes the proof.
% \end{beweis}

In the following, interpreting AL6ab8 geometrically, we construct two parabolas and two tangents to them such that two given points are superposed by folding across these tangents. Finding these tangents resp. folds $l_1$ and $l_2$ as in Figure 2 is equivalent to finding the intersection points of two (origami) cubics. 

By the Bézout Theorem there are nine (projective) complex intersection points of two cubics. The equation~(\ref{cubic}) reveals that, in our situation, the two origami cubic curves will have two intersection points at infinity, so we get generically seven affine intersection points. This yields an equation of degree 7, cf.\ the following section.

The Galois group of this equation is $S_7$ in the generic case, but we noticed that, for suitable values for the points and lines, smaller group such as $A_7$ and $PSL_3\ff_2$ occur as well (see Figures 4 and 5 for concrete crease patterns for these groups). This observation led to the natural question 
\emph{Which subgroups of $S_7$ are realizable by 2-fold axioms?} Even stronger: \emph{Is every septic equation solvable by 2-fold axioms?}

\section{Generic Septic Equations}
\label{generic}
% One might think that finding a solution to the problem  ``Is every septic equation solvable by 2-fold axiom?'' is purely computational: One has to 
% There are rational polynomials realising $PSL_{3}\ff_{2}$ for instance, but calculating with degree-7-polynomials in several unknowns exceeds magmas possibilities, this time. So we had to think of some simplifications.
We want to find suitable values for the given points and lines, such that the arising origami cubic curves intersect in a point with the ``right'' minimal polynomial, i.e. minimal polynomial with the wanted Galois group, for instance. We have seen that it suffices in many cases to fix one of the two parabolas.
 In the axiom AL6ab8 we drop the generality and specify one half of the data. Let  $n$ be the line with the equation $y=-1$. Let $R=(0,1)$ and $S=(0,0)$, cf. Figure 2. So  we fold the origin across the tangents of the parabola with the equation $y = \frac{1}{4}x^2$. If $(X,Y)$ is an intersection point of the two arising origami cubic curves and $W\!:=\frac{X}{Y}$ (this is, up to sign, just the slope of the fold line across which $(0,0)$ is reflected to $(X,Y)$),
the following equation of degree $7$ is satisfied by $W$:
\begin{equation}\label{deg7pol}
\begin{gathered}
 W^7 + \big(\frac{3}{2}e + \frac{1}{2}ft_0 + t_0 - \frac{1}{2}t_1\big)W^6+\\ + \big(\frac{3}{4}e^2 + \frac{1}{2}eft_0 + et_0 -\frac{1}{2}et_1 + \frac{1}{4}f^2 - \frac{1}{4}ft_2 + f - \frac{1}{2}t_2\big)W^5+\\
+ \big(\frac{1}{8}e^3 + \frac{1}{8}e^2ft_0 + \frac{1}{4}e^2t_0 - \frac{1}{8}e^2t_1 + \frac{1}{8}ef^2 - \frac{1}{8}eft_2 +\frac{1}{2}ef - \frac{1}{4}et_2 + \frac{1}{2}e + \frac{1}{8}f^3t_0+ \\+ \frac{3}{4}f^2t_0 - \frac{1}{8}f^2t_3 +\frac{3}{2}ft_0 - \frac{1}{2}ft_3 - \frac{1}{2}t_3\big)W^4+\\ + \big(\frac{3}{4}e^2 + \frac{1}{2}eft_0 - \frac{1}{2}et_1 + \frac{1}{4}f^2 - \frac{1}{4}ft_2\big)W^3+\\ + \big(\frac{1}{4}e^3 + \frac{1}{4}e^2ft_0 + \frac{1}{4}e^2t_0 -\frac{1}{4}e^2t_1 + \frac{1}{4}ef^2 - \frac{1}{4}eft_2+ \\ + \frac{1}{2}ef -\frac{1}{4}et_2 + \frac{1}{4}f^3t_0 +
        \frac{3}{4}f^2t_0 - \frac{1}{4}f^2t_3 - \frac{1}{2}ft_3\big)W^2+\\ + \frac{1}{8}e^3 +\frac{1}{8}e^2ft_0 -
        \frac{1}{8}e^2t_1 + \frac{1}{8}ef^2 - \frac{1}{8}eft_2 + \frac{1}{8}f^3t_0 - \frac{1}{8}f^2t_3 = 0.
\end{gathered}
\end{equation}
We see that the coefficient at $W^1$ is zero (which can of course always be achieved for a general equation of degree $7$ by substituting $W^{-1}$ for $W$ and applying a linear transformation).
The question is therefore whether the remaining six coefficients can take arbitrary values $s_1,\ldots,s_6$ for suitable choices of $t_0,t_1,t_2,t_3,e,f$. As the resulting system of 
equations is linear in $t_0,\ldots,t_3$, we can assign arbitrary values to four of the coefficients (say, the coefficients at $W^6,\ldots, W^3$).\\
More precisely, this is achieved by setting
\begin{align*}
&t_0\!:=\frac{2s_1e + s_2f - \frac{3}{2}e^2 - \frac{1}{2}f^2-s_4(f+2)}{ef+2e},\\
&t_1\!:=3e + ft_0 + 2t_0 -2s_1,\\
&t_2\!:=\frac{4s_1e - 3e^2 + f^2  + 4f -4s_2}{2+f},\\
&t_3\!:=\frac{-2s_1e^2 + 4s_2e + e^3 + 4e + f^3t_0 + 6f^2t_0 + 12ft_0 -8s_3}{4f+4+f^2}.\end{align*}
The polynomial we obtain in this way from Formula (\ref{deg7pol}) is
\begin{equation}
\label{deg7pol2}
W^7+s_1W^6+s_2W^5+s_3W^4+s_4W^3+C_1W^2+C_2 = 0,
\end{equation}
where 
\begin{align*}
C_1:=\frac{1}{4(e f^2 + 4 e f + 4 e)}\big(&-2s_1 e^3 f - 4s_1 e^3 - 2 s_1 e f^3 - 12 s_1 e f^2 - s_2 e^2 f^2 + 2 s_2 e^2 f + 8 s_2 e^2 -  s_2 f^4 -\\&- 6 s_2 f^3 + 8 s_3 e f^2 +16 s_3 e f +  s_4 e^2 f^2 + 4s_4 e^2 f + 4s_4 e^2 + s_4 f^4 + 8 s_4 f^3 +\\
        &+12 s_4 f^2 + \frac{1}{2} e^4 f +  e^4 + e^2 f^3 + 4e^2 f^2 - 8 e^2 f +
        \frac{1}{2} f^5 + 3 f^4\big),\end{align*}
\begin{align*}
C_2:=\frac{1}{4(e f^3 + 6 e f^2 + 12 e f + 8 e)}\big(&-2 s_1 e^3 f^2 - 4s_1 e^3 f - 2 s_1 e f^4 - 8 s_1 e f^3 -  s_2 e^2 f^3 +
        4s_2 e^2 f -\\&-  s_2 f^5 - 4s_2 f^4 + 4s_3 e f^3 + 8 s_3 e f^2 +  s_4 e^2 f^3
        + 6 s_4 e^2 f^2 +\\&+ 12 s_4 e^2 f + 8 s_4 e^2 +  s_4 f^5 + 6 s_4 f^4 +
        8 s_4 f^3 + \frac{1}{2} e^4 f^2 - 2 e^4 +\\ &  +e^2 f^4 + 2 e^2 f^3 -6 e^2 f^2 + \frac{1}{2} f^6 + 2 f^5\big).
\end{align*}
Now we are ready to show the main result.

\begin{Satz}
 A generic equation of degree $7$ can be solved by 2-fold origami.
\end{Satz}
\begin{beweis}

If we set $f=0$ in equation (\ref{deg7pol2}), then we obtain
\begin{equation}
\label{deg7pol_reduced}
W^7 + s_1W^6 + s_2W^5 + s_3W^4 + s_4W^3 + \frac{1}{16}(e^3 - 4e^2s_1 + 8es_2 + 8es_4)W^2 - \frac{1}{16}e^3 + \frac{1}{4}es_4=0.
\end{equation}
Replacing $W$ by $W^{-1}$, we obtain a septic equation with vanishing coefficient at $W^6$. After multiplying $W$ with an appropriate factor and dividing by the leading coefficient, we
even get a monic septic polynomial of the form $W^7+a_1W^5+a_2W^4+a_3W^3+a_4W^2+a_5W+a_5$, where the $a_i$ are rational functions in $s_1,\ldots,s_4$ and $e$.\\
We investigate whether for suitable choices of $s_1,\ldots,s_4$ and $e$ \emph{any} equation of the form
 \begin{equation}\label{general_septic} W^7+a_1W^5+a_2W^4+a_3W^3+a_4W^2+a_5W+a_5,\end{equation} with real-valued coefficients $a_1,\ldots,a_5$,
can be obtained. This leads to a system of polynomial equations in the variables $s_1,\ldots,s_4$ and $e$ over the function field $\qq(a_1,\ldots,a_5)$. Gröbner basis methods show that
the system can be solved by $e$ satisfying the equation
 \begin{align*}
p_e(a_1,\ldots,a_5)\!:= e^8&  + \\
\frac{e^6}{a_5^3}&\big(48a_1^2a_2a_5^2 + 40a_1a_2^4a_5 - 176a_1a_2^2a_5^2 - 56a_1a_4a_5^2 + 112a_1a_5^3 + 4a_2^7-\\
 &- 40a_2^5a_5 - 40a_2^3a_4a_5 + 128a_2^3a_5^2 + 136a_2a_4a_5^2 - 128a_2a_5^3\big)+\\
 \frac{e^4}{a_5^4}&\big(368a_1^4a_2^2a_5^2 + 448a_1^4a_5^3 - 64a_1^3a_2^5a_5 - 192a_1^3a_2^3a_5^2 -
        1248a_1^3a_2a_4a_5^2 - 1280a_1^3a_2a_5^3 +\\
	&+ 448a_1^2a_2^4a_4a_5 - 96a_1^2a_2^4a_5^2 + 1888a_1^2a_2^2a_4a_5^2 + 128a_1^2a_2^2a_5^3 + 896a_1^2a_4^2a_5^2 - 896a_1^2a_4a_5^3 +\\
	&+ 1792a_1^2a_5^4 + 128a_1a_2^5a_4a_5 - 64a_1a_2^5a_5^2 - 1184a_1a_2^3a_4^2a_5 -
        544a_1a_2^3a_4a_5^2 + 512a_1a_2^3a_5^3 - \\
	&-320a_1a_2a_4^2a_5^2 - 640a_1a_2a_4a_5^3 - 1024a_1a_2a_5^4 - 32a_2^6a_4^2 + 64a_2^6a_4a_5 - 16a_2^6a_5^2 + \\
	&+ 352a_2^4a_4^2a_5 - 736a_2^4a_4a_5^2 + 192a_2^4a_5^3 + 800a_2^2a_4^3a_5 - 1600a_2^2a_4^2a_5^2 + 2816a_2^2a_4a_5^3 -\\
        &- 768a_2^2a_5^4 - 896a_4^3a_5^2 + 3584a_4^2a_5^3 - 3584a_4a_5^4 + 1024a_5^5\big)+\\
  \frac{e^2}{a_5^5}&\big(256a_1^7a_5^3 + 1024a_1^6a_2a_5^3 - 1024a_1^5a_2^2a_4a_5^2 + 1536a_1^5a_2^2a_5^3 - 3584a_1^5a_4a_5^3 -\\
	&- 2048a_1^4a_2^3a_4a_5^2 + 1024a_1^4a_2^3a_5^3 + 4608a_1^4a_2a_4^2a_5^2 -
        8704a_1^4a_2a_4a_5^3 + 512a_1^3a_2^4a_4^2a_5 - \\
	&- 1024a_1^3a_2^4a_4a_5^2 + 256a_1^3a_2^4a_5^3 + 12800a_1^3a_2^2a_4^2a_5^2 - 6656a_1^3a_2^2a_4a_5^3 - 3584a_1^3a_4^3a_5^2 + \\
	&+ 14336a_1^3a_4^2a_5^3 - 3072a_1^2a_2^3a_4^3a_5 + 6144a_1^2a_2^3a_4^2a_5^2 -
        1536a_1^2a_2^3a_4a_5^3 - 27648a_1^2a_2a_4^3a_5^2 +\\
	&+ 15360a_1^2a_2a_4^2a_5^3 + 7168a_1a_2^2a_4^4a_5 - 12288a_1a_2^2a_4^3a_5^2 + 3072a_1a_2^2a_4^2a_5^3 + 14336a_1a_4^4a_5^2 - \\
	& - 14336a_1a_4^3a_5^3 + 64a_2^5a_4^4 - 768a_2^3a_4^4a_5 - 4608a_2a_4^5a_5 + 11264a_2a_4^4a_5^2 - 2048a_2a_4^3a_5^3\big)+\\
 \frac{1}{a_5^5}&\big(-1024a_1^3a_2^3a_4^4 + 6144a_1^2a_2^2a_4^5 - 12288a_1a_2a_4^6 + 8192a_4^7\big)=0
 \end{align*}
and $s_1,\ldots,s_4$ lying in the field extension of $\qq(a_1,\ldots,a_5)$ generated by $e$.

 But obviously $e^2$  is a root of  a quartic polynomial. As quadratic and quartic polynomials can be solved by 1-fold origami,
 $e$ is an origami-constructible number – and so are $s_1,\ldots,s_4$. Therefore, by substituting $t_0,\ldots,t_3$ and then $a,b,c,d$ as described above, all the values for our 2-fold step
 are constructible numbers. If we can, in addition, choose them as real numbers – for which it is sufficient that $e$ is real – then we can solve the generic septic equation (\ref{general_septic}) by 2-fold origami.

While the above polynomial $p_e(a_1,\ldots,a_5)$ of degree $8$ may of course have no real roots for certain choices of $a_1,\ldots,a_5$, we will show that there is always a polynomial
\[W^7+b_1W^5+b_2W^4+b_3W^3+b_4W^2+b_5W+b_5\] generating the same field extension as the analogous polynomial in $a_1,\ldots,a_5$, such that $p_e(b_1,\ldots,b_5)$ has a real root.

Firstly, observe that $p_0(a_1,\ldots,a_5) = -1024\cdot a_4^4a_5^{-5}\cdot (a_1a_2-2a_4)^3$ and $\lim\limits_{e\to+\infty}p_e= +\infty$. If we can enforce $a_5(a_1a_2-2a_4)>0$,
 then $p$ will change its sign somewhere between $0$ and $+\infty$ and therefore have a real root. Now for $w\in \rr$ a root of
$W^7+a_1W^5+a_2W^4+a_3W^3+a_4W^2+a_5W+a_5$ and $\lambda\in \qq$, we can bring the minimal polynomial of $w+\frac{\lambda}{w}$ into the form $W^7+b_1W^5+b_2W^4+b_3W^3+b_4W^2+b_5W+b_5$
via linear transformations. The term $b_5(b_1b_2-2b_4)$ is a rational function  in the $a_i$ and $\lambda$; as we are only interested in the sign of this expression, we can multiply it by arbitrary squares and thus obtain a square-free polynomial  $F$ in $a_1,\ldots, a_5$ and $\lambda$.

 Viewing $F$ as a polynomial in $\lambda$ over $\qq(a_1,\ldots, a_5)$, we observe that $F$ splits as $F(\lambda)=F_1(\lambda)\cdot F_2(\lambda)$ with  $F_1$, $F_2$ polynomials in $\lambda$ of degree 5 and 7 respectively. But $F_1$ and $F_2$ will both have a real root, and generically these roots will not coincide; this means that the expression $b_5(b_1b_2-2b_4)$
will change its sign at some point, so if we choose $\lambda\in \qq$ in a suitable interval, $b_5(b_1b_2-2b_4)$ will be positive, and $p_e(b_1,\ldots,b_5)$ will have a real root.
But this means that we can construct $w+\frac{\lambda}{w}$, and therefore $w$ as well, with 2-fold origami, so every real root of a generic septic equation is constructible by 2-fold origami. 
\end{beweis}
\begin{Bemerkung}
 Note that our ``generic'' form can be obtained without loss of generality, if we view the coefficients as
 transcendentals; however, for certain specializations, like polynomials of the form $W^7-A$ this is not possible by linear transformations. We will deal with equations $W^7-A = 0$ in \ref{seventhroot}.\\
Also, throughout the proof, we deal with rational functions in certain coefficients; of course, for a bad choice of the coefficients, these might not be well-defined due to vanishing denominators. The term ``generic'' polynomial should always be understood in the sense that the denominators have to behave well. 
\end{Bemerkung}

\section{Solvable groups}
We showed above that a generic equation of degree $7$ is solvable by 2-fold origami, but there are some important cases which seem not to be included in the generic result,
like 2-folding of seventh roots. We deal with this separately and show more generally that every solvable $\{2,3,5,7\}$-extension of $\qq$ is solvable by 2-fold origami.
\subsection{Angle septisection}
If you are an origami artist you have quite often to create some difficult marks to proceed. Usually these are some divisions of a segment, like third parts.
It can occur that you need a third part of an angle\footnote{By the way, the possibility of angle trisection is one of the advantages of 1-fold origami over euclidean constructions.}. 
Robert Lang found an exact angle quintisection with 2-fold origami, which is impossible by 1-fold origami, and \cite{AL} and \cite{Ni} put this result on a more general basis. As far as we know an exact angle septisection for a general angle has not been given by means of $k$-fold origami for $k<5$. 
Robert Lang did find an approximate solution \cite{lang2010}, though, and used it for the construction of his famous scorpion.

Let $\varphi\in (0,2\pi)$ be an angle, $A=2\cos(\varphi)$ and $x=2\cos(\varphi/7)$. Then one easily verifies with de~Moivre's formula that $x^7-7x^5+14x^3-7x-A=0$.\\
If we can solve this equation for arbitrary $A\in (-2,2)$, then we can septisect an arbitrary angle. The following theorem states that we can do this with 2-fold origami.
\begin{Satz}
Septisection of arbitrary angles $\varphi\in(0,2\pi)$ is possible with 2-fold origami. 
\end{Satz}
\begin{beweis}
We take the polynomial from equation (\ref{deg7pol_reduced}), replace $W$ with ${W}^{-1}$ (so the polynomial has vanishing coefficient at $W^6$ instead of $W^1$), and multiply $W$ with a constant factor
in order to let the constant and the linear coefficient take the same value. Denote the resulting polynomial by $h_1(W)$. Then we treat $W^7-7W^5+14W^3-7W-A$ in the same way (that is, multiply $W$ with factor $\frac{A}{7}$) and denote the result by $h_2(W)$. Now compare the coefficients of
$h_1$ and $h_2$. The arising system of equations over $\qq(A)$ is solved by $s_2=0=s_4$ and  
\begin{align*}
s_1:=\;&\frac{4302592 (-28 + A^2) (-196 + 14 A^2 + 3 A^4) e + 
 196 A^4 (5488 + 560 A^2 + A^4) e^3 - A^6 (28 + 3 A^2) e^5}{153664 A^2 (21952 - 784 A^2 - 252 A^4 + A^6)},\\
s_3:=\;&\frac{-3764768 (-112 + A^4) e - 98 A^2 (784 + 280 A^2 + A^4) e^3 + A^6 e^5}{5488 (21952 - 784 A^2 - 252 A^4 + A^6)},
\end{align*}
where $e$ fulfills $e^6 - \frac{38416}{A^2} e^4 + \frac{-7529536 A^4 + 210827008 A^2 - 843308032}{A^6} e^2 - \frac{210827008}{A^2} = 0$.
As all the other unknown coefficients $a,b,c,d$ of our initial point and line setting can be expressed as rational functions in these, we are done if we can construct $e$ as a real number; but the above sextic polynomial in $e$ can be solved by solving cubic and quadratic equations, i.\;\!e.\ by 1-fold origami. It remains to be seen whether $e$ can be chosen as a real number. 
As $p(x) = x^6 - \frac{38416}{A^2} x^4 + \frac{-7529536 A^4 + 210827008 A^2 - 843308032}{A^6} x^2 - \frac{210827008}{A^2}$ is negative at $0$ and $\lim\limits_{x\to+\infty} p(x) = +\infty$, such a real number $e$ exists, indeed.
% Regard $\qq(A,e)\mid \qq(A)$ as a function field extension (with a transcendental $A$). The number of real solutions can then only change at a ramification point, but the only ramification point inside the interval $(-2,2)$ is at $A=0$. As only even powers of $A$ appear in our equation, the number of real solutions for values of $A$ left and right of the ramification point must be the same; it suffices to verify for one specialization (say $A=1$) that this number is 2. This concludes the proof. 
\end{beweis}

\subsection{Folding seventh roots}\label{seventhroot}
We try to specialize all intermediate coefficients of the polynomial in equation (\ref{deg7pol2}) to zero. This corresponds to constructing seventh roots.
%For many values of the remaining constant coefficient, the corresponding systems of equations become solvable over small fields. We give the values that lead to the field $\qq(\sqrt[7]{2})$. 
%\marginpar{Werte}\\
So we compare coefficients of the polynomial in (\ref{deg7pol2}) with those of the polynomial $W^7+s$, where $s$ is any positive real number. 
This leads to two equations in $e$ and $f$ over the field $\qq(s)$. This system of equations has a solution in the function field  defined by\\
$f^{10}t^2 + 2f^{10}t + f^{10} + 24f^9t^2 + 24f^9t + 252f^8t^2 - 84f^8t + 1536f^7t^2 - 1264f^7t + 6048f^6t^2 - 1008f^6t + 16128f^5t^2 + 5376f^5t + 29568f^4t^2 + 3584f^4t + 36864f^3t^2 - 6144f^3t + 
         29952f^2t^2 + 14336ft^2 + 3072t^2=0,$ where $t\!:=s^2$. This defines a rational function field $\qq(f,t)$ over $\qq(t)$, and therefore we can find a parameter $w$ such that $\qq(w)=\qq(f,t)$ and express $t$ as a rational
function in it; computer calculation yields 
$t=2^{10}\frac{w^7}{(w + 7)^7(w+1)^2(w+3)}$ for a suitable parameter $w$.

Remember that we want to solve $X^7+\sqrt{t}=0$. Multiply $X$ with a factor $\sqrt{\frac{2w}{w+7}}$, we can transform this to $X^7+\sqrt{T}=0$, where 
$T = \frac{8}{(w+1)^2(w+3)}$. Note that the square root that is introduced in this transformation does not lead to any problems, as square roots are of course constructible by 1-fold origami.

But now we can specialize $T$ to an arbitrary positive value; $w$ will then be the (w.l.o.g.\ real) root of a cubic equation, and we can solve this equation with 1-fold origami. Now $e$ and $f$ lie in the field generated by $w$ and $\sqrt{t}$, which is at most a  quadratic extension of $\qq(w)$. As we can w.\;\!l.\;\!o.\;\!g.\ multiply $T$ with positive rational 7th powers, the field $\qq(w,\sqrt{T})$ can even be enforced to be real because  for $T>0$ small enough, for the equation $8 = T(w+1)^2(w+3)$ will always have a \emph{positive} solution $w$, and therefore $t$ will be positive with $T$ as well. So the construction is completed. 
%
%
%Remember that we want to solve $X^7+\sqrt{t}$. Stretching $X$ with a factor $\sqrt{\frac{2}{7(w +1/7)}}$, we can transform this to $X^7+\sqrt{\tilde{t}}$, where 
%$\tilde{t} = 8/3\frac{w^3}{(w+1)^2(w+1/3)}$.\\
%Note that the square root that is introduced in this transformation does not lead to any problems, as square roots are of course constructible by (standard one-fold) origami.\\
%But now we can specialize $\tilde{t}$ to an arbitrary positive value; $w$ will then be the (w.l.o.g.\ real) root of a cubic equation, and we can solve this equation with one-fold 
%origami; as $e$ and $f$ lie in the field generated by $w$ and $\sqrt{t}$, which is at most a (real)\footnote{For the reality, note that we can w.l.o.g.\ multiply $\tilde{t}$ with 
%positive rational 7-th powers; for $\tilde{t}>0$ small enough, the equation $8/3w^3 = \tilde{t}(w+1)^2(w+1/3)$ will have a positive solution, 
%and therefore $t$ will be positive with $\tilde{t}$ as well.} quadratic extension of $\qq(w)$, we are done. 
%
 
Together with angle septisection shown above, this result leads to the following
\begin{Satz}
Let $K\mid \qq$ be a finite solvable Galois extension of degree $2^a\cdot 3^b\cdot 5^c\cdot 7^d$ with $a,b,c,d\in \nn_0$. Then $K$ is solvable by 2-fold origami.
\end{Satz}
\begin{beweis}

Galois theory says that the extension $K\mid \qq$ can be solved by repeatedly taking (square, cubic, fifth and seventh) roots. 
Now taking the $n$-th root of any complex number can be achieved by taking the real $n$-th root of its absolute value, combined with angle $n$-section.\\
Square roots and cubic roots can be taken by 1-fold origami. Nishimura \cite{Ni} and Lang \cite{lang5} showed that in particular fifth roots and quintisection can be taken with 2-fold origami. This leaves $n=7$, and we showed above how to 
septisect arbitrary angles and take seventh roots of reals.
\end{beweis}

\section{{Crease patterns for nonsolvable transitive groups in $S_7$}}
In the previous section we showed that every polynomial whose Galois group is a solvable subgroup of $S_7$ can be solved by 2-fold origami. 
Now we turn to nonsolvable transitive groups in $S_7$. These are $S_7$, $A_7$ and $PSL_3\ff_2\cong PSL_2\ff_7$, cf. \cite[p.\ 60, Table 2.1]{dixon}. With the methods of Section \ref{generic}, one could give many explicit constructions
for each of these groups; however these constructions would in general be quite lengthy and involved as they require for instance the folding of solutions of quartic equations.\\
We give explicit examples of folds with very nice initial coordinates that lead to Galois groups $A_7$ and $PSL_3\ff_2$ (the generic case $S_7$ is left out as almost all
folds with axiom AL6ab8 lead to this Galois group).

First, we want to give a realisation of $A_{7}$ by specializing the axiom AL6ab8. 
\begin{figure}[!ht]\centering
\includegraphics[height=8cm]{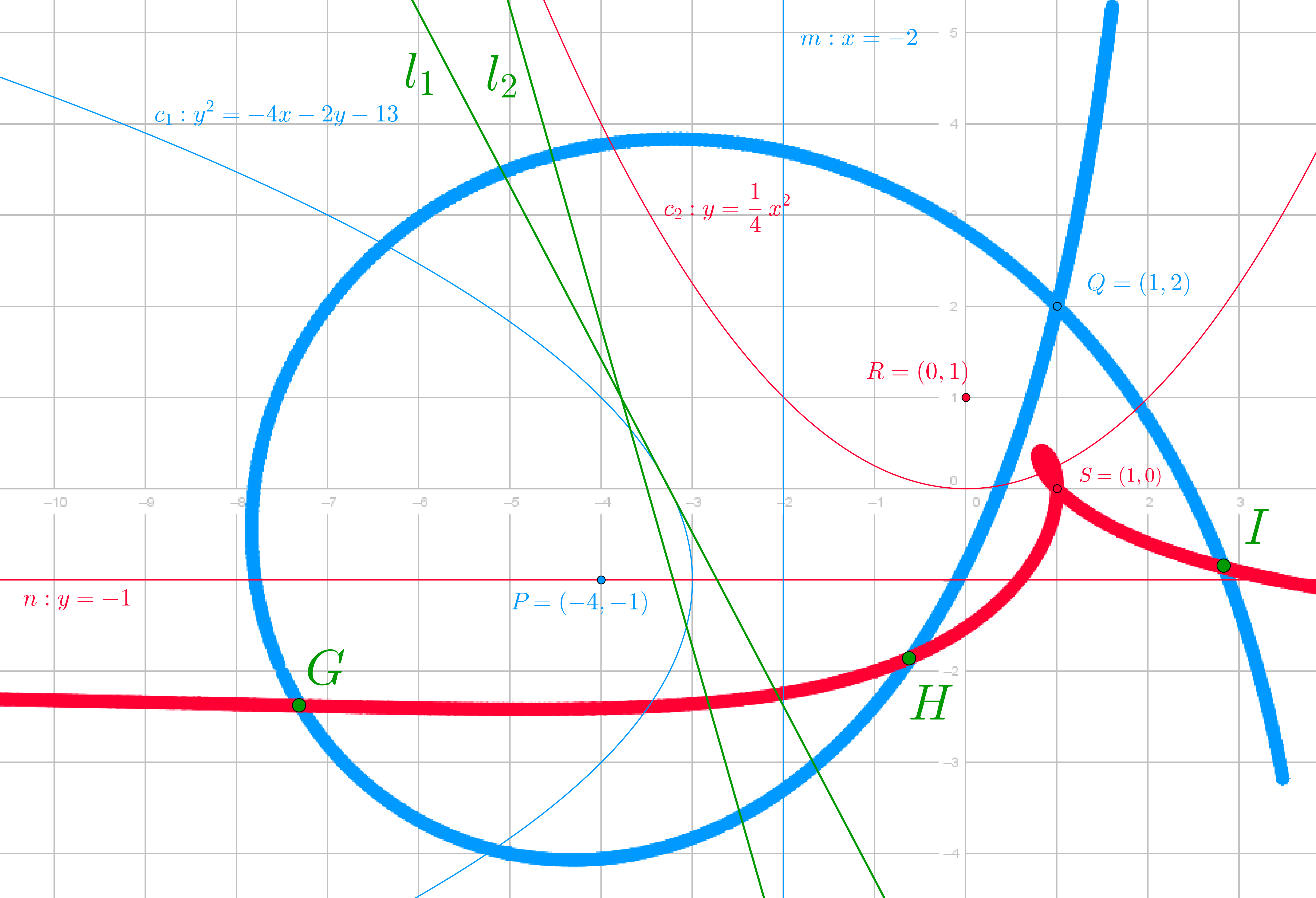}
 \caption{Crease pattern for $A_7$ by AL6ab8. $G,H,I$ are the intersection points of the two bold cubics in red and blue. 
 The green foldlines $l_1$ and $l_2$ arise by folding $Q$ resp. $S$ on $G$.}
\end{figure}
We put \[m: x= -2 ,\; P=(-4,-1),\; Q=(1,2)\] for the first parabola set, cf.\ Figure 4. Furthermore we set \[n: y = -1 ,\; R=(0,1),\; S=(1,0)\] for the second parabola set. 
Putting these numbers into the equations we dealt with above, we get a polynomial $h$ of degree 7, describing the intersection points of the two cubics, such that $\operatorname{Gal}(h\mid\qq) \cong A_{7}$.
More precisely, the slope of the foldline $l_2$ is a root of the polynomial $y^7 + y^6 - 8y^5 + 3y^4 + y^3 - 3y^2 + 2y - 1$. The discriminant of this polynomial is equal to 
$2^8\cdot 31^2\cdot 157^2$, so it is a square and the Galois group must be contained in $A_7$. In fact, equality holds, as one verifies with a computer algebra program such as Magma.\\
Note that this polynomial has exactly three real roots, corresponding to the three intersection points of our cubics in the affine plane. The slope of the line $l_2$ in Figure 4 is the real root of approximate value $-3.49$.

Now, let us describe how to construct $PSL_{3}\ff_{2}$ by AL6ab8. As depicted in Figure 5, set \[m: y = \frac{1}{2}x-1,\; P=(-\frac{16}{5}, -\frac{12}{5}),\; Q=(-3,-3); \quad n: y = -2,\; R=(0,0),\; S=(1,-1).\]
\begin{figure}[!ht]\centering
\includegraphics[height=7cm]{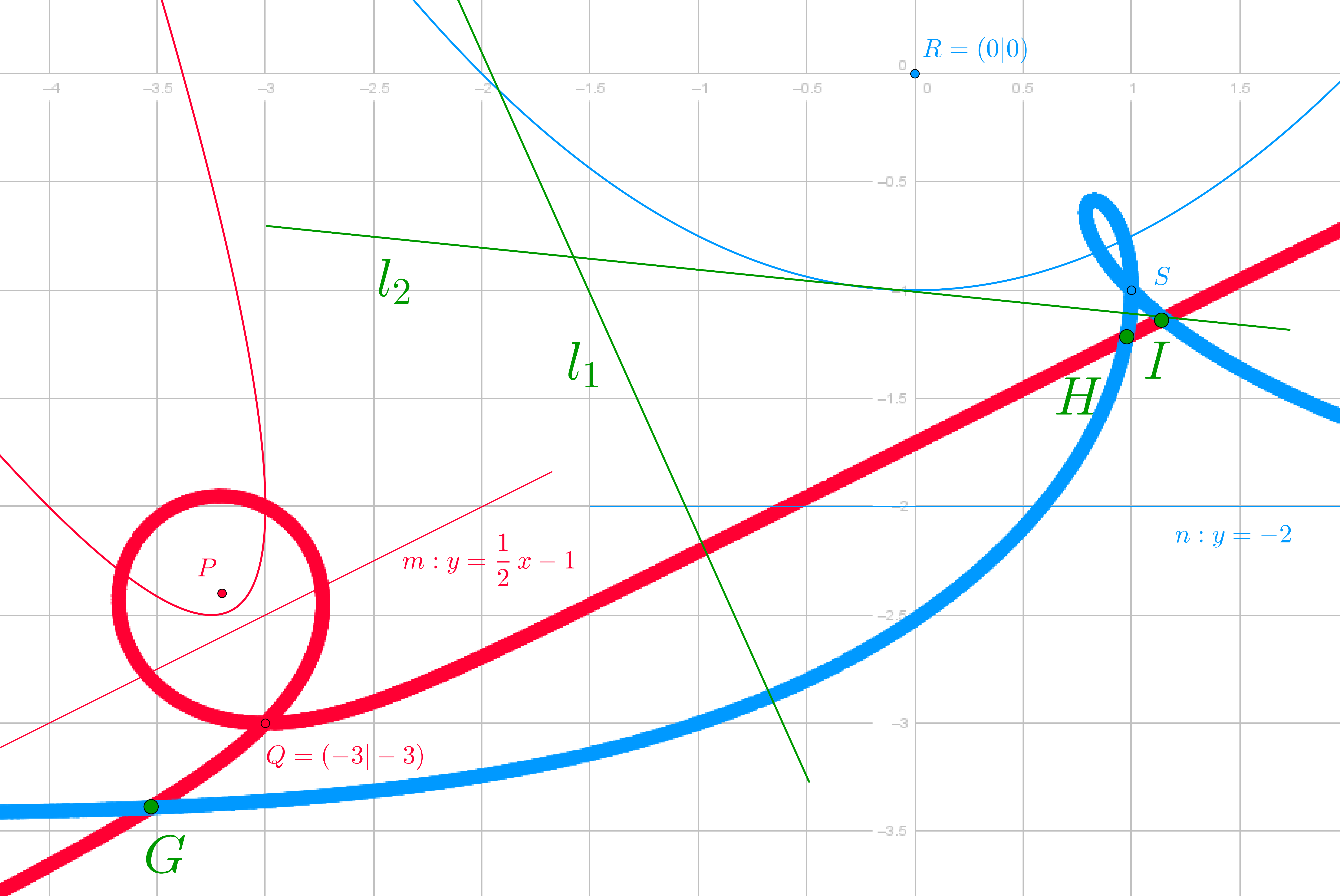}
 \caption{Crease pattern for $PSL_3\ff_2$ by AL6ab8. The green foldlines $l_1$ and $l_2$ arise by folding $Q$ resp. $S$ on $H$.}
\end{figure}
Again the two cubics intersect in three real points; the slope of fold line $l_2$ fulfills the equation \[y^7 + 3y^6 - 3y^4 + 5y^3 + y^2 - 10y - 1=0,\] 
whose Galois group surprisingly turns out to be $PSL_{3}\ff_{2}$. It is notable that this polynomial is very simple and the number field generated by one of its roots has very small
discriminant, namely $2^6\cdot 383^2$.

\end{document}